\documentclass[12pt]{amsart}
\usepackage{amssymb}
\usepackage{amsthm}
\usepackage{enumerate}

\renewcommand{\phi}{\varphi}

\newcommand{\kt}{K_\Theta}

\newcommand{\vep}{\varepsilon}

\newcommand{\rs}{\tz\setminus \rho(\Theta)}

\newcommand{\dd}{\mathbb{D}}
\newcommand{\tz}{\mathbb{T}}

\renewcommand{\Re}{{\rm Re}\,}

\newtheorem{Thm}{Theorem}[section]
\newtheorem{Lem}[Thm]{Lemma}

\newtheorem{Cor}[Thm]{Corollary}

\def\beginpf{\noindent {\bf Proof.} \ }

\small\normalsize

\title[The Feichtinger conjecture for reproducing kernels]
{The Feichtinger conjecture for reproducing \\
kernels in model subspaces}
\author{Anton Baranov \& Konstantin Dyakonov}

\address{Department of Mathematics and Mechanics, St. Petersburg State University,
28, Universitetskii pr., St. Petersburg, 198504, Russia}
\email{anton.d.baranov@gmail.com}

\address{ICREA and Universitat de Barcelona, Departament de Matem\`atica
Aplicada i An\`alisi, Gran Via 585, E-08007 Barcelona, Spain}
\email{dyakonov@mat.ub.es}
\keywords{Hardy space, inner function, model subspace, reproducing kernel,
Riesz basis, Bessel sequence, Feichtinger conjecture}
\subjclass[2000]{30D50, 30D55, 30E05, 46E22.}
\thanks{The first author was supported in part by grants MK 5027.2008.1
and NSH 2409.2008.1 (Russia). The second author was supported in
part by grants MTM2008-05561-C02-01 and MTM2007-30904-E from El
Ministerio de Ciencia e Innovaci\'on (Spain).}

\begin{document}
\begin{abstract}
We obtain two results concerning the Feichtinger
conjecture for systems of normalized reproducing kernels in the model
subspace $K_\Theta = H^2\ominus \Theta H^2$ of the Hardy space $H^2$,
where $\Theta$ is an inner function. First, we verify the Feichtinger conjecture
for the kernels $ \tilde k_{\lambda_n} = k_{\lambda_n}/\|k_{\lambda_n}\|$
under the assumption that $\sup_n |\Theta(\lambda_n)|<1$. Secondly, we prove
the Feichtinger conjecture in the case where $\Theta$ is a one-component
inner function, meaning that the set $\{z:|\Theta(z)|<\varepsilon\}$ is
connected for some $\varepsilon\in(0,1)$.
\end{abstract}

\maketitle

\section{Introduction}

A sequence of unit vectors $\{h_n\}$ in a separable Hilbert space $\mathcal H$
is said to be a {\it Bessel sequence} if, for some constant $C>0$
and every $h\in\mathcal H$,
\begin{equation}
\label{0}
\sum_n\left|\langle h,h_n\rangle_{\mathcal H}\right|^2 \le C\|h\|^2_{\mathcal H}.
\end{equation}
Further, a sequence $\{h_n\}\subset\mathcal H$ is called a {\it Riesz basic sequence} if
it is a Riesz basis in its span, or equivalently, if
there exists a constant $A>0$ such that
$$
A^{-1}\sum_n|c_n|^2 \le \Big\|\sum_n c_n h_n\Big\|^2_{\mathcal H}
\le A\sum_n|c_n|^2
$$
for any finite sequence $\{c_n\}\subset\mathbb C$.

The Feichtinger conjecture states that {\it every Bessel sequence
splits into finitely many Riesz basic sequences.}

The Feichtinger conjecture is a problem of high current interest.
As recently shown in \cite{Cas1, Cas2, Cas3}, it is
equivalent to the famous Kadison--Singer conjecture in $C^*$-algebras
and to some other important open problems in analysis.

As an attempt to better understand the heart of the problem
(or to find a counterexample), one may look at the Feichtinger conjecture
for special systems in function spaces. A natural class of examples
is given by systems of (normalized) reproducing kernels
in Hilbert spaces of analytic functions. Let $\mathcal H$ be a
reproducing kernel Hilbert space of analytic functions on some domain $D$,
and let $k_\lambda$ denote the kernel function corresponding to a point $\lambda\in D$,
so that  $f(\lambda)=\langle f,k_\lambda\rangle_{\mathcal H}$ for each $f\in\mathcal H$.
The Bessel property (\ref{0}) for a system of normalized reproducing kernels
$k_{\lambda_n}/\|k_{\lambda_n}\|_{\mathcal H}$ is equivalent to the
Carleson-type embedding $\mathcal H\subset L^2(\mu)$ for the discrete measure
$$
\mu = \sum_n \|k_{\lambda_n}\|^{-2}_{\mathcal H}\, \delta_{\lambda_n},
$$
where $\delta_{\lambda_n}$ is the unit point mass at $\lambda_n$.
Carleson measures are well understood for many
classical spaces (e.g., for Hardy, Bergman and Bargmann--Fock spaces),
and then the validity of the Feichtinger conjecture follows
from various known results about sampling and interpolation in
these spaces.
\bigskip

\section{Main results}

In this paper we consider a class of reproducing kernel Hilbert spaces
where the problem is still open. Let $H^2$ denote the Hardy space of the
unit disk $\mathbb{D}$, equipped with the standard norm
$\|\cdot\|_2 = \|\cdot\|_{L^2(m)}$, where $m$ is normalized
Lebesgue measure on the unit circle $\mathbb{T}$.
The reproducing kernels are then the usual Cauchy kernels
$(1-\overline \lambda z)^{-1}$, and
the Feichtinger conjecture is true. Indeed,
if $\mu = \sum_n (1-|\lambda_n|)\, \delta_{\lambda_n}$
is a Carleson measure, then $\Lambda = \{\lambda_n\}$ is a finite union of
interpolating sequences, and each of these corresponds to a Riesz basic sequence
of Cauchy kernels. However, the problem becomes nontrivial for
{\it model} (or {\it star-invariant}) {\it subspaces} of $H^2$,
that is, for subspaces of the form
$$
K_\Theta=H^2\ominus \Theta H^2,
$$
where $\Theta$ is an inner function on the disk.
These subspaces play a distinguished role in operator theory (see,
e.g., \cite{Nikshift, Nik}) and in operator-related complex analysis.

The reproducing kernel for $K_\Theta$
corresponding to a point $\lambda\in\mathbb{D}$ is given by
$$
k_\lambda(z)=
\frac{1-\overline{\Theta(\lambda)}\Theta(z)}{1-\overline \lambda z}.
$$
Since functions in $K_\Theta$ have more analyticity than general $H^2$
functions, there may
exist reproducing kernels at boundary points (see Section 3
for details). We shall denote the normalized kernel
$k_{\lambda}/\|k_{\lambda}\|$ by $\tilde k_\lambda$.

The geometry of reproducing kernels in model subspaces
seems to be highly nontrivial. No complete description is known
either for Bessel sequences or for Riesz bases, which makes
the Feichtinger conjecture for these spaces a difficult problem.
The problem of describing the Carleson measures for $K_\Theta$ was
posed by Cohn in 1982 and is still open, even though a number of partial results
are available \cite{al2, bar05, con1, DSpb, DAmer, DJFA, vt}. One special case
where the Carleson measures for $K_\Theta$ have been completely characterized is
that of a {\it one-component} inner function $\Theta$. Here, by saying that $\Theta$
is one-component we mean that the set $\{z\in \dd: |\Theta(z)|<\vep\}$
is connected for some $\vep\in (0,1)$.

The Riesz bases and Riesz basic sequences of reproducing kernels in $K_\Theta$
were described by Hruscev, Nikol'skii and Pavlov (see \cite{hnp} and
also \cite[Part D, Chapter 4]{Nik}) under the additional hypothesis that
the sequence $\Lambda=\{\lambda_n\}$ be close to the {\it spectrum} $\rho(\Theta)$
of $\Theta$, in the sense that
\begin{equation}
\label{1}
\sup_n |\Theta(\lambda_n)|<1.
\end{equation}
Recall that the spectrum $\rho(\Theta)$ of an inner function $\Theta$
is, by definition, the smallest closed subset of $\mathbb D\cup\mathbb T$ containing
the zeros of $\Theta$ and the support of the associated singular measure.
Equivalently, $\rho(\Theta)$ is the set of all points
$\zeta\in\mathbb D\cup\mathbb T$
such that $\liminf\limits_{z\to\zeta,\, z\in\mathbb{D}}|\Theta(z)|=0$.

To summarize the Hruscev--Nikol'skii--Pavlov results, suppose that $\Lambda=\{\lambda_n\}$
obeys (\ref{1}) and let $B_\Lambda$ be the Blaschke product with zero
sequence $\Lambda$. Then $\{\tilde k_{\lambda_n}\}$ is a Riesz basis in
$K_\Theta$ if and only if $\Lambda$ satisfies the {\it Carleson condition}
$$
({\rm C})\qquad\qquad
\delta(B_\Lambda) :=  \inf_n \prod_{k:\,k\ne n} \bigg|\frac{\lambda_n-\lambda_k}
{1-\overline\lambda_k\lambda_n}\bigg| >0\qquad\qquad\text{\,\,\,\,\,\,\,\,\,\,\,\,}
$$
(i.e., $\Lambda$ is an {\it interpolating sequence})
and the Toeplitz operator $T_{\Theta \overline B_{\Lambda}}$
is invertible (this can be further rephrased by invoking the Devinatz--Widom
invertibility criterion; see, e.g., \cite[p.\,234]{hnp}). Similarly,
$\{\tilde k_{\lambda_n}\}$ is a Riesz basic sequence
if and only if $\Lambda$ satisfies (C) and $T_{\Theta \overline B_{\Lambda}}$
is left-invertible; the latter condition is equivalent
to ${\rm dist} (\Theta \overline B_{\Lambda}, H^\infty) <1.$

We are now in a position to state our main results. The first of these
applies to a general inner function $\Theta$, provided that
the Hruscev--Nikol'skii--Pavlov condition (\ref{1}) is fulfilled. The underlying
Hilbert space is, of course, always taken to be $K_\Theta$.

\begin{Thm}
Every Bessel sequence of normalized reproducing kernels $\{\tilde k_{\lambda_n}\}$
with property (\ref{1}) splits into finitely many Riesz basic sequences.
\end{Thm}

The second theorem treats the case of a one-component inner function.

\begin{Thm}
Assume that $\Theta$ is a one-component inner function. Then every Bessel
sequence of normalized reproducing kernels $\{\tilde k_{\lambda_n}\}$
splits into finitely many Riesz basic sequences.
\end{Thm}

We shall arrive at Theorem 2.2 by combining Theorem 2.1 with some stability results
(essentially due to Cohn) for Riesz bases of reproducing kernels.

It should be noted that our results carry over, via a unitary transform,
to model subspaces of the Hardy space $H^2(\mathbb{C}_+)$ in the upper half-plane
$\mathbb{C}_+$. In this case, the premier example of a one-component inner function
is given by $\Theta_a(z):=\exp(iaz)$, $a>0$, and $K_{\Theta_a}$ essentially reduces
to the Paley--Wiener space $PW_a$. Precisely speaking, if $PW_a$ stands
for the set of entire functions of exponential type at most $a$ that are
square integrable on $\mathbb R$, then
$$K_{\Theta_a} = PW_a\cap H^2(\mathbb{C}_+) = e^{iaz/2} PW_{a/2},$$
and we conclude that the Feichtinger conjecture holds for
reproducing kernels in $PW_a$. This statement follows also
from Theorem 2.1, since the sequence $\{\lambda_n\}$
generates a Riesz sequence of normalized kernels in $PW_a$
if and only if the same is true for the sequence  $\{\lambda_n + i\}$
(which satisfies (\ref{1}) when $\lambda_n \in
\mathbb{C_+} \cup \mathbb{R}$).

Recall that the Fourier transform identifies $PW_a$ with $L^2(-a,a)$
and converts the reproducing kernels $k_\lambda$ into the exponentials
$e_\lambda(t):=\exp(i\lambda t)$ on the interval $(-a,a)$. Thus, we have

\begin{Cor}
The Feichtinger conjecture is true for any
Bessel sequence of normalized exponentials
$\{e_{\lambda_n}/\|e_{\lambda_n}\|\}$ in $L^2(-a,a)$.
\end{Cor}

Quite recently, Yu. S. Belov, T. Y. Mengestie and K. Seip \cite{bms} have 
described the Bessel sequences and proved the Feichtinger conjecture for a class 
of model subspaces in $H^2(\mathbb{C}_+)$ generated by Blaschke products with
very sparse zeros. Such inner functions are, of course, never one-component.
Thus, our current results and those of \cite{bms} complement each other.

The rest of the paper is organized as follows. In the next section we prove
Theorem 2.1. Section 4 contains some preliminaries on Clark bases and a discussion of
their stability. Finally, in Section 5, these results are applied to prove Theorem 2.2.

In what follows, the letter $C$ will stand for a positive constant, not necessarily
the same in different places. We write  $A\asymp B$ to mean that $C^{-1}B\le A\le CB$
for some constant $C>0$.
\bigskip

\section{Proof of Theorem 2.1}

It is known that the Bessel property for a system of normalized 
reproducing kernels $\{\tilde k_{\lambda_n}\}$ always implies that 
the measure $\sum_n (1-|\lambda_n|)\, \delta_{\lambda_n}$ is a Carleson
measure (see \cite[Part D, Lemma 4.4.2]{Nik} or \cite[Lemma 4.2]{bar05a}).
This is especially easy to see when (\ref{1}) is satisfied.
Indeed, observe that
$$
\|k_\lambda\|_2^2=k_\lambda(\lambda)=\frac{1-|\Theta(\lambda)|^2}{1-|\lambda|^2},
\qquad\lambda\in\mathbb D.
$$
Then (\ref{1}) implies that $\|k_{\lambda_n}\|_2^2 \asymp (1-|\lambda_n|^2)^{-1}$
and $|k_{\lambda_m}(\lambda_n)| \asymp |1-\overline\lambda_m\lambda_n|^{-1}$.
Applying the Bessel inequality
\begin{equation}\label{eqn:bessel}
\sum_n|f(\lambda_n)|^2\left\|k_{\lambda_n}\right\|_2^{-2}\le C\|f\|_2^2,\qquad f\in K_\Theta,
\end{equation}
to $f=k_{\lambda_m}$, we deduce that
$$
\sup_m \sum_n \frac{(1-|\lambda_m|^2)(1-|\lambda_n|^2)}
{|1-\overline\lambda_m\lambda_n|^2} <\infty.
$$
The latter condition means (cf. \cite[p.\,151]{Nikshift}) that
$\sum_n (1-|\lambda_n|)\, \delta_{\lambda_n}$ is a Carleson measure, or equivalently
(see \cite[p.\,158]{Nikshift}), that $\Lambda$ is a finite union of interpolating 
sequences. Clearly, no generality will be lost in assuming that $\Lambda$ is a single
interpolating sequence.

Recall that, for an interpolating sequence $\Lambda=\{\lambda_n\}$, the
{\it constant of interpolation} $c(\Lambda)$ is defined as the smallest constant $c$
with the following property: whenever $\{a_n\}\subset\mathbb C$ and
$\sup_n|a_n|\le1$, there exists a function $f\in H^\infty$ with $\|f\|_\infty\le c$
that solves the interpolation problem $f(\lambda_n)=a_n$ ($n=1,2,\dots$).
Recall also that the interpolation constant $c(\Lambda)$ tends to 1 from above,
as the Carleson constant $\delta = \delta(B_\Lambda)$ in (C) tends to 1
from below. More explicitly, an estimate due to Earl \cite{earl} reads
\begin{equation}\label{eqn:earlest}
1\le c(\Lambda)\le\varphi\left(\delta(B_\Lambda)\right),
\end{equation}
where
$$
\varphi(\delta):=\frac{2-\delta^2 + 2(1-\delta^2)^{1/2}}{\delta^2}.
$$

Our aim is to split $\Lambda$, a given interpolating sequence, into finitely many (say $N$)
subsequences $\Lambda_j=\left\{\lambda_n^j\right\}_{n=1}^\infty$ such that
\begin{equation}\label{eqn:dist}
{\rm dist}\,(\Theta, B_j H^\infty)<1,\qquad j=1,\dots,N,
\end{equation}
where $B_j$ is the Blaschke product with zeros $\Lambda_j$.
This done, we shall readily conclude (by the discussion in Section 2 above) that the kernels
corresponding to the points of $\Lambda_j$ form a Riesz basic sequence, for each $j$.
We shall thus arrive at the sought-after partition.

To see how \eqref{eqn:dist} can be achieved, note that
$$
{\rm dist}\,(\Theta, B_j H^\infty)=\inf \big\{ \|f\|_\infty: f\in H^\infty,
f|_{\Lambda_j}=\Theta|_{\Lambda_j}\big\} \le \gamma c_j,
$$
where $\gamma:=\sup_n|\Theta(\lambda_n)|<1$ and $c_j:=c(\Lambda_j)$ is the
interpolation constant for the subsequence $\Lambda_j$. Therefore, it suffices
to make sure that
\begin{equation}\label{eqn:gammacj}
c_j<1/\gamma\quad\text{\rm for}\quad j=1,\dots,N.
\end{equation}

Using Earl's estimate \eqref{eqn:earlest} with $\Lambda_j$ in place of $\Lambda$, 
we see that \eqref{eqn:gammacj} becomes true provided that the Carleson
constants $\delta_j:=\delta\left(B_{\Lambda_j}\right)$ get close enough
to 1, so as to ensure $\varphi(\delta_j)<1/\gamma$. This in turn can be arranged via
(a corollary of) Mills' factorization lemma, as stated in \cite[Chapter X, Corollary 1.6]{ga}.
The lemma tells us that every interpolating Blaschke product $B$ can be factored as $B=B_1B_2$
so that $\delta(B_k)\ge (\delta(B))^{1/2}$, $k=1,2$. Applying the same procedure to
each of the factors that arise, and iterating this as many times as we need, we arrive
at subproducts $B_j$ whose $\delta_j$'s are as close to 1 as required.
\qed
\bigskip

\section{Clark bases and stability}

There may exist Riesz bases of reproducing kernels
which do not satisfy (\ref{1}).
An important example is given by a Clark basis
corresponding to points of $\mathbb{T}$.

Let us begin by recalling that, by Ahern and Clark's results \cite{ac},
we have $k_\zeta\in K_\Theta$ for a point $\zeta\in \tz$ if and only if
\begin{equation}
\label{3}
|\Theta'(\zeta)|=\sum\limits_n \frac{1-|z_n|^2}{|\zeta-z_n |^2}+
2\int\limits_\mathbb{T}\frac{d\nu(\tau)}{|\zeta - \tau|^2}<\infty.
\end{equation}
Here $z_n$ are the zeros of $\Theta$ and $\nu$ is the associated
singular measure.

Now we turn to Clark's construction
of orthogonal bases of reproducing kernels \cite{cl}.
For each $\alpha\in\mathbb{T}$, the function $(\alpha+\Theta)
/(\alpha-\Theta)$ has positive real part in $\mathbb{D}$.
Hence, there exists a finite (singular) positive measure
$\sigma_\alpha$ on $\tz$ such that
$$
\Re \frac{\alpha+\Theta(z)}{\alpha-\Theta(z)}=
\int\limits_\mathbb{T} \frac{1-|z|^2}{|\tau - z|^2}\,
d \sigma_\alpha(\tau), \qquad z\in\mathbb{D}.
$$
Clark's theorem states that if $\sigma_\alpha$ is purely atomic,
i.e., if $\sigma_\alpha = \sum_n a_n\, \delta_{\tau_n}$,
then the system $\{k_{\tau_n}\}$ is an orthogonal basis in ${K^2_{\Theta}}$;
in particular, $k_{\tau_n}\in {K^2_{\Theta}}$.
Note that all measures $\sigma_\alpha$ are
purely atomic when the boundary spectrum $\rho(\Theta)\cap\mathbb{T}$
is at most countable.
\smallskip

For a one-component inner function $\Theta$,
it was shown by Aleksandrov \cite{al2} that the set
$\rho(\Theta)\cap \tz$ has zero Lebesgue measure,
and moreover, $\sigma_\alpha(\rho(\Theta) \cap \tz)=0$
for every Clark measure $\sigma_\alpha$.
Since $\tz \setminus \rho(\Theta)$ is a countable union of arcs
where $\Theta$ is analytic, each Clark measure
is atomic and generates an orthogonal (Clark) basis of
reproducing kernels $k_{\tau_n}$; here $\{\tau_n\}$ is an
enumeration of the level set
$$
\{\tau\in\tz\setminus \rho(\Theta): \,\Theta(\tau)=\alpha\}.
$$
Note that $\Theta$ has an increasing smooth branch of
the argument on each arc of $\tz\setminus \rho(\Theta)$. Furthermore,
the change of argument of $\Theta$ between
two neighboring points $\tau_n$ and $\tau_{n+1}$ equals $2\pi$. That is,
$$
\int_{(\tau_n,\tau_{n+1})}|\Theta'(\tau)|dm(\tau)
=\frac1{2\pi i}\int_{(\tau_n,\tau_{n+1})}\frac{\Theta'(\tau)}{\Theta(\tau)}d\tau
=1,
$$
where $(\tau_n,\tau_{n+1})$ is the corresponding arc on $\tz$
and $m$ is {\it normalized} Lebesgue measure on $\tz$.

We shall need the following estimate due to Aleksandrov
\cite{al2}: if $\Theta$ is a one-component inner function, then
there exists a positive constant $C = C(\Theta)$ such that
$$
\Bigg|\frac{1-\Theta(z)\overline{\Theta(\tau)}}
{(\tau -z)\cdot\Theta'(\tau)}\Bigg|\le C
$$
for all $z\in \dd$ and $\tau\in \tz \setminus \rho(\Theta)$.
Since $\liminf\limits_{z\to\zeta,\, z\in\mathbb{D}}|\Theta(z)|=0$
for $\zeta\in \rho(\Theta)$, this implies
\begin{equation}
\label{3a}
    |\Theta'(\tau)|^{-1} \le C {\rm dist}\, (\tau, \rho(\Theta)).
\end{equation}
We see, in particular, that $\int_J |\Theta'(\tau)|\,dm(\tau) =\infty$
for every connected component $J$ of $\tz \setminus \rho(\Theta)$,
unless $\Theta$ is a finite Blaschke product.

An interesting result about stability of Clark bases
for one-component functions was obtained by Cohn \cite[Theorem 3]{cohn86}.
It says that $\{\tilde k_{\lambda_n}\}$ will be a Riesz basis in $K_\Theta$
provided that $\{\lambda_n\}$ is close to the support $\{\tau_n\}$
of a Clark basis, in the sense that the variation of $\Theta$ between
$\tau_n$ and $\lambda_n$ is small.

\begin{Thm}
Suppose $\Theta$ is a one-component inner function, $\{k_{\tau_n}\}$ is
a Clark basis in $K_\Theta$, and $\lambda_n \in \mathbb{D} \cup \tz$.
Then there is an $\varepsilon=\varepsilon(\Theta)>0$
with the property that $\{\tilde k_{\lambda_n}\}$
is a Riesz basis in $\kt$ whenever
there exist paths $(\tau_n,\lambda_n)$ connecting
$\tau_n$ and $\lambda_n$ for which
$$
\sup_n \int_{(\tau_n,\lambda_n) } |\Theta'(\tau)|\, |d\tau| < \vep.
$$
\end{Thm}

This will be an important ingredient in our proof of Theorem 2.2.
In fact, we shall use the following slight modification of Theorem 4.1,
which is a particular case of \cite[Corollary 1.3]{bar05a}.

\begin{Thm}
Let $\Theta$ be a one-component inner function.
Then there exists an $\varepsilon= \varepsilon(\Theta) >0$
making the following statement true: if $\{k_{\tau_n}\}$ is a Clark basis
in ${K_{\Theta}}$ and if $\lambda_n$ are points of $\dd\cup\tz$ with
$$
|\lambda_n-\tau_n|<\varepsilon |\Theta'(\tau_n)|^{-1},
$$
then $\{\tilde k_{\lambda_n}\}$ is a Riesz basis in $\kt$.
\end{Thm}

It should be emphasized that these results depend heavily on
properties of {\it one-component} inner functions and fail in the general
case; see \cite{bar05a} for counterexamples and for an extension that holds
for generic inner functions.
\bigskip

\section{Proof of Theorem 2.2}

Throughout this section, $\Theta$ is a one-component inner function and
$\{\tilde k_{\lambda_n}:\lambda_n\in\Lambda\}$ is a Bessel
sequence in $K_\Theta$. The idea of the proof is to split $\Lambda$ into
two sequences, one of which is contained in the sublevel set
$\{|\Theta|<\delta\}$ for some $\delta<1$, while the other
is a small perturbation of certain Clark measures' supports.
We shall use a special system of arcs and Carleson squares
that was introduced in \cite{bar08}, where compactness
and Schatten class properties of the embeddings $K_\Theta \subset L^2(\mu)$
were studied.

Given a large positive integer $N$, we define the (countable)
sets $T_l$ with $l=1, \dots, N$ by
$$
T_l=\{\tau_m^l\}=\big\{\tau\in\rs: \,\Theta(\tau)=
e^{2\pi i l/N}\big\}.
$$
Each $T_l$ corresponds to a Clark basis.
We also consider the set $\{\zeta_n\}=
\bigcup\limits_{l=1}^N T_l$. Then we have a partition of
$\tz\setminus \rho(\Theta)$ into arcs $J_n$ with mutually
disjoint interiors,
whose endpoints are in the set $\{\zeta_n\}$ (we always assume that
$\zeta_n$ is the first endpoint of $J_n$ when moving clockwise)
and which satisfy
\begin{equation}
\label{2}
\int_{J_n}|\Theta'(\tau)|\,dm(\tau)=1/N.
\end{equation}

\begin{Lem}
If $N$ is sufficiently large, then
$$
|\Theta'(\zeta)|\asymp|\Theta'(\zeta_n)|,\qquad\zeta\in J_n,
$$
where the constant involved is numerical. In particular,
\begin{equation}
\label{4}
|J_n|\asymp N^{-1}|\Theta'(\zeta_n)|^{-1},
\end{equation}
where $|J_n|$ denotes the length of the arc $J_n$.
\end{Lem}

\beginpf
By (\ref{2}), $|J_n|=2\pi N^{-1} |\Theta'(\tau)|^{-1}$ for
some $\tau\in J_n$. It follows then from (\ref{3a}) that, for
$N$ suitably large,
$$
|J_n|\le C {\rm dist}\, (\tau, \rho(\Theta))/N < {\rm dist}\,
(\tau, \rho(\Theta))/10,
$$
so that the length of the arc is much smaller than the distance
to the spectrum. A trivial estimate based on formula (\ref{3})
gives us $|\Theta'(\zeta)|\asymp|\Theta'(\tau)|$
for $\zeta \in J_n$ (with an absolute constant),
and (\ref{4}) follows.
\qed
\bigskip

With each arc $J_n$ we associate the {\it Carleson square}
$$
S_n:=\{r\zeta:\,\zeta \in J_n,\,1-|J_n|/2\pi \le r\le 1\},
$$
and we put
$$
G:=\dd\setminus \bigcup\limits_{n} S_n.
$$
Further, for a fixed $l$, we write $J_m^l$ ($m=1,2,\dots$) to enumerate those arcs among
the $J_n$'s whose first endpoint (when moving clockwise) is in $T_l$, and we denote the
corresponding Carleson squares by $S_m^l$.

Obviously, ${\rm diam}\, S_n\asymp |J_n|\asymp N^{-1}|\Theta'(\zeta_n)|^{-1}$.
Theorem 4.2 therefore implies the following

\begin{Cor} If $N$ is sufficiently large and $l$ is any fixed index
in $\{1,\dots,N\}$, then $\{\tilde k_{\lambda_m}:m=1,2,\dots\}$ is a Riesz basis
in $\kt$ whenever $\{\lambda_m\}$ is a sequence with $\lambda_m\in S_m^l$.
\end{Cor}

From now on, we fix some value of $N=N(\Theta)$ that makes the conclusions of
Lemma 5.1 and Corollary 5.2 true. Our next step is to prove the following lemma.

\begin{Lem}
There exists a $\delta=\delta(\Theta)$, $0<\delta<1$,
such that $|\Theta(z)|<\delta$ for all $z\in G$.
\end{Lem}

\beginpf
Let $\zeta \in \mathbb{T}\setminus \rho(\Theta)$ and let
$z\in \mathbb{D}$ be a point with
\begin{equation}\label{eqn:distzz}
|z-\zeta| < {\rm dist}\, (\zeta, \rho(\Theta))/2.
\end{equation}
An elementary estimate then yields
\begin{equation}
\label{4a}
\log|\Theta(z)| \le -C(1-|z|)|\Theta'(\zeta)|
\end{equation}
with an absolute constant $C>0$ (see, e.g., the proof of Theorem 4.9
in \cite{bar05}).

Consider the boundary $\partial G$ of the domain $G$, and
note that $\partial G\cap \tz = \rho(\Theta)\cap \tz$.
Our plan is to show that there exists a $\delta\in (0,1)$ such that
$|\Theta(z)|\le \delta$ for each $z\in \partial G\cap \dd$. Since
$\partial G$ is a rectifiable Jordan curve, $|\Theta(z)|\le 1$ in $G$, and
$\partial G\cap\tz$ is of zero Lebesgue measure, the desired conclusion
that $|\Theta(z)|<\delta$ for all $z\in G$ is then guaranteed by the maximum
principle.

Now suppose $z\in \partial G\cap \dd$. There are two possibilities:
either, for some $n$, $z=(1-|J_n|/(2\pi))\, \zeta$ with $\zeta \in J_n$
(i.e., $z$ lies on the interior side of some square $S_n$)
or there are two adjacent squares
$S_m$ and $S_n$ with $|J_m|\le |J_n|$ such that
$z=r\zeta$, where 
$$
1-|J_n|/(2\pi) \le r\le 1-|J_m|/(2\pi)
$$
and $\zeta$ is the common endpoint of $J_m$ and $J_n$.
Note that, by Lemma 5.1, $|J_m|\asymp |J_n|
\asymp N^{-1}|\Theta'(\zeta)|^{-1}$. Thus, in any case,
$$
1-|z|\asymp |J_n|\asymp N^{-1}|\Theta'(\zeta)|^{-1}, \qquad
\zeta = z/|z|,
$$
with some absolute constants, while \eqref{eqn:distzz} holds true.

Consequently, (\ref{4a}) implies $\log|\Theta(z)|\le -C/N$, and so
$|\Theta(z)|\le e^{-C/N}=:\delta$ on $\partial G\cap\dd$.
\qed
\bigskip

To treat the points in  $\Lambda\cap \bigcup\limits_m S_m$, we are going 
to show that, for a Bessel sequence $\{\tilde k_{\lambda_n}\}$,
the number of $\lambda_n$'s in each of the squares is uniformly bounded. 
But first we need to establish a certain \lq\lq monotonicity property" of 
the norm $\|k_{r\zeta}\|_2$ as a function of $r\in (0,1)$. 

\begin{Lem}
There is an absolute constant $C>0$ such that 
\begin{equation}\label{monoton}
\|k_{r\zeta}\|^2_2\le C|\Theta'(\zeta)|\left(=C\|k_\zeta\|^2_2\right)
\end{equation}
whenever 
$\zeta\in\tz\setminus\rho(\Theta)$ and $0<r<1$. 
\end{Lem}

\beginpf
One easily checks (\ref{monoton}) when $\Theta$ is a single Blaschke 
factor $b_{\eta}(z):=
\frac{|\eta|}{\eta}\cdot \frac{\eta-z}{1-\overline{\eta}z}$. If $\Theta =
\prod_j b_{z_j}$ is a Blaschke product, then (\ref{monoton}) is due 
to the fact that
$$
\|k_{r\zeta}\|_2^2 = \frac{1-|\Theta(r\zeta)|^2}{1-r^2} \le
\sum\limits_j \frac{1-|b_{z_j}(r\zeta)|^2}{1-r^2} \le C\sum_j|b'_{z_j}(\zeta)|
=C|\Theta'(\zeta)|.
$$
Finally, the general case follows from Frostman's theorem on 
approximation of an arbitrary inner function by Blaschke products.
\qed
\bigskip

An extension (and a different proof) of Lemma 5.4 can be found 
in \cite[Corollary 4.7]{bar05}.

\begin{Lem}
There is a constant $M$ such that each set $\Lambda\cap S_m$ consists
of at most $M$ points.
\end{Lem}

\beginpf
Let $w_m$ be the midpoint of the interior side (i.e., the smaller circular arc) of 
$\partial S_m$. We have then $1-|w_m|=|J_m|/(2\pi)$ and $|\Theta(w_m)|\le\delta<1$, 
where $\delta$ is the same as in Lemma 5.3. Applying the Bessel inequality 
\eqref{eqn:bessel} with $f=k_{w_m}$ and noting that
 $$\|k_{w_m}\|_2^2 \le (1-|w_m|)^{-1}=2\pi|J_m|^{-1},$$
we obtain
\begin{equation}\label{eqn:bessker}
\sum_n\bigg|\frac{1-\overline{\Theta(w_m)}
\Theta(\lambda_n)}{1-\overline w_m\lambda_n}\bigg|^2\|k_{\lambda_n}\|_2^{-2}
\le C|J_m|^{-1}.
\end{equation}
Observe, in addition, that $|1-\overline{\Theta(w_m)}\Theta(\lambda_n)|\ge1-\delta$,
while for $\lambda_n \in S_m$ we also have
$$
|1-\overline w_m\lambda_n|\le C|J_m|\qquad\text{\rm and}\qquad
\|k_{\lambda_n}\|_2^2\le C|J_m|^{-1}.
$$
Here, the last inequality relies on Lemma 5.4, combined with the fact that 
$|\Theta'(\zeta)|\asymp|J_m|^{-1}$ for $\zeta\in J_m$; the constant involved 
depends only on $\Theta$. 

Therefore, multiplying \eqref{eqn:bessker} by $|J_m|$ and restricting the
summation to the indices $n$ with $\lambda_n\in S_m$ yields
$$
\#\{n: \lambda_n\in S_m\} \le C\sum_{\lambda_n\in S_m}
\frac{(1-\delta)^2|J_m|} {\|k_{\lambda_n}\|_2^2\,|1-\overline w_m\lambda_n|^2}
\le\text{\rm const}.
$$
We conclude that
$\#\{n: \lambda_n\in S_m\}$ is bounded by a constant independent of $m$.
\qed
\bigskip

Now we complete the proof of Theorem 2.2.
To deal with the family $\{\tilde k_{\lambda_n}: \lambda_n\in G\}$
we may apply Theorem 2.1, since in this case
$|\Theta(\lambda_n)|\le \delta<1$ by Lemma 5.3.

Then we split $\Lambda \cap \bigcup\limits_m S_m $
into $N$ sets $\{\lambda_n^l\}$, $l=1, \dots, N$,
such that $\{\lambda_n^l\}\subset \bigcup\limits_j S_j^l$,
where $S_j^l$ are the Carleson squares corresponding to the $l$th
Clark basis $\{k_\tau:\tau\in T_l\}$.
Finally, for each $l$ we write $\{\lambda_n^l\}$ as a union of
at most $M$ sets $\{\lambda_n^{l,m}\}$
with the property that each Carleson square $S_j^l$ contains at most
one point from $\{\lambda_n^{l,m}\}$. Now Corollary 5.2 shows
that the normalized reproducing kernels corresponding to each set
$\{\lambda_n^{l,m}\}$ form a Riesz basic sequence.
\qed
\bigskip
\\
{\bf Acknowledgements.}
We are grateful to Kristian Seip, Tesfa Mengestie and Yurii Belov for 
their interest in this work and for supplying us with a prepublication 
version of their paper \cite{bms}. The second author is also indebted to
Joaquim Ortega-Cerd\`a for a stimulating discussion; in particular,
it was his idea to look at the Feichtinger conjecture for $\kt$ in the 
case where $\sup_n|\Theta(\lambda_n)|<1$.

\end{document}